\documentclass{amsart}
\usepackage{amssymb}

\newtheorem{theorem}{Theorem}[section]
\newtheorem{proposition}{Proposition}[section]
\newtheorem{corollary}{Corollary}[section]
\newtheorem{lemma}{Lemma}[section]

\newtheorem*{KC}{Kaplansky's Conjecture}
\newtheorem*{QTC}{Quasitrace Conjecture}

\newtheorem*{claimone}{\sl Claim 1}
\newtheorem*{claimtwo}{\sl Claim 2}

\newtheorem*{partcase}{\sl Particular Case}

\theoremstyle{definition}
\newtheorem*{definition}{Definition}

\newtheorem*{comment}{Comment}

\newtheorem*{notations}{Notations}
\newtheorem{remark}{Remark}[section]
\newtheorem*{remarks}{Remarks \theremark}
\newtheorem{example}{Example}[section]
\newtheorem*{notation}{Notation}

\newcommand{\dontwrite}[1]{}

\newcommand{\twone}{\mbox{$\text{\rm II}_1$}}
\newcommand{\onefin}{\mbox{$\text{\rm I}_{\text{\rm fin}}$}}
\newcommand{\typeone}[1]{\mbox{$\text{\rm I}_{\text{$#1$}}$}}

\begin{document}
\title[Abelian Self-Commutators in Finite Factors]{Abelian Self-Commutators
in
Finite Factors}
\author{Gabriel Nagy}
\address{Department of Mathematics, Kansas State university, Manhattan KS 66506, U.S.A.}
\email{nagy@math.ksu.edu}
\begin{abstract}
An {\em abelian self-commutator\/} in a C*-algebra
$\mathcal{A}$ is an element of the form $A=X^*X-XX^*$, with
$X\in\mathcal{A}$, such that $X^*X$ and $XX^*$ commute.
It is shown that, given a finite AW*-factor $\mathcal{A}$, there
exists another finite AW*-factor $\mathcal{M}$ of same type as $\mathcal{A}$,
that contains $\mathcal{A}$ as an
AW*-subfactor, such that any self-adjoint element $X\in\mathcal{M}$
of quasitrace zero is an abelian self-commutator in $\mathcal{M}$.
\end{abstract}
\maketitle

\section*{Introduction}

According to the Murray-von Neumann classification, finite
von Neumann factors are either of type I${}_{\text{fin}}$, or of type
\twone. For the non-expert, the easiest way to understand this classification is
by accepting the famous result of Murray and von Neumann (see \cite{MvN})
which states that {\em every finite von Neumann factor $\mathcal{M}$
 posesses a  unique state-trace $\tau_{\mathcal{M}}$}.
Upon accepting this result, the type of $\mathcal{M}$ is decided by so-called
{\em dimension range\/}:
$\mathcal{D}_{\mathcal{M}}=\big\{\tau_{\mathcal{M}}(P)\,:\,P\text{ projection in
}\mathcal{M}\big\}$ as follows. If $\mathcal{D}_{\mathcal{M}}$ is finite, then
$\mathcal{M}$ is of type I${}_{\text{fin}}$ (more explictly, in this case
$\mathcal{D}_{\mathcal{M}}=\big\{\frac kn\,:\,k=0,1,\dots,n\big\}$ for some
$n\in\mathbb{N}$, and $\mathcal{M}\simeq \text{Mat}_n(\mathbb{C})$ -- the 
algebra of $n\times n$ matrices). If $\mathcal{D}_{\mathcal{M}}$ is infinite, then
$\mathcal{M}$ is of type \twone, and in fact one has $\mathcal{D}_{\mathcal{M}}=[0,1]$.
From this point of view, the factors of type \twone\ are the ones
that are interesting, one reason being the fact that, although all factors of type
\twone\ have
the same dimension range, there are uncountably many non-isomorphic ones
(by a celebrated result of Connes).

In this paper we deal with the problem of characterizing
the self-adjoint elements
of trace zero, in terms of simpler ones. We wish to carry this investigation
in a ``Hilbert-space-free'' framework, so instead of von Neumann factors, we
are going to work within the category of AW*-algebras.
Such objects were introduced in the 1950's by Kaplansky (\cite{Kap}) in an
attempt to formalize the theory of von Neumann algebras without any use of
pre-duals. Recall that
A unital C*-algebra $\mathcal{A}$ is called an {\em AW*-algebra}, if for every non-empty
set $\mathcal{X}\subset\mathcal{A}$, the left anihilator
set
$\mathbf{L}(\mathcal{X})=\big\{A\in\mathcal{A}\,:\,AX=0,\,\,\,
\forall\,X\in\mathcal{X}\big\}$
is the principal right ideal generated by a projection $P\in\mathcal{A}$, that is,
 $\mathbf{L}(\mathcal{X})=\mathcal{A}P$.

Much of the theory -- based on the geometry of projections -- works for
AW*-algebras exactly as in the von Neumann case, and
one can classify the finite AW*-factors into the types
I${}_{\text{fin}}$ and \twone, exactly as above, but using the
following alternative result: {\em any finite AW*-factor $\mathcal{A}$ posesses
a unique normalized quasitrace $q_{\mathcal{A}}$}. Recall that a
{\em quasitrace\/} on a C*-algebra $\mathfrak{A}$ is a map $q:\mathfrak{A}\to\mathbb{C}$ with the
following properties:
\begin{itemize}
\item[(i)] if $A,B\in\mathfrak{A}$ are self-adjoint, then
$q(A+iB)=q(A)+iq(B)$;
\item[(ii)] $q(AA^*)=q(A^*A)\geq 0$, $\forall\,A\in\mathfrak{A}$;
\item[(iii)] $q$ is linear on all abelian C*-subalgebras of $\mathfrak{A}$,
\item[(iv)] there is a map $q_2: \text{Mat}_2(\mathcal{A})\to\mathbb{C}$ with
properties (i)-(iii), such that
$$
q_2\left(\left[
\begin{array}{cc}
A& 0\\ 0 &0
\end{array}
\right]\right)=q(A),
\,\,\,\forall\,A\in\mathfrak{A}.
$$
\end{itemize}
(The condition that $q$ is {\em normalized\/} means that $q(I)=1$.)

With this terminology, the dimension range of a finite AW*-factor is the set
$\mathcal{D}_{\mathcal{A}}=\big\{q_{\mathcal{A}}(P)\,:\,P\text{ projection in
}\mathcal{A}\big\}$, and the classification into the two types
is eaxctly as above. As in the case of von Neumann factors, one can show that
the AW*-factors of type I${}_{\text{fin}}$ are again the matrix algebras
$\text{Mat}_n(\mathbb{C})$, $n\in\mathbb{N}$. The type \twone\ case however is
still mysterious. In fact, a longstanding problem in the theory of AW*-algebras
is the following:
\begin{KC}
Every AW*-factor of type \twone\ is a von Neumann factor.
\end{KC}
An equivalent formulation states that: {\em if $\mathcal{A}$ is an AW*-factor of type
\twone, then the quasitrace $q_{\mathcal{A}}$ is linear\/}
(so it is in fact a trace). It is well known (see \cite{Ha} for example)
that Kaplansky's Conjecture implies:
\begin{QTC}
Quasitraces (on arbitrary C*-algebras) are traces.
\end{QTC}
A remarkable result of Haagerup (\cite{Ha}) states that quasitraces on
{\em exact\/} C*-algebras are traces, so
if $\mathcal{A}$
is an AW*-factor
of type \twone, generated (as an AW*-algebra) by an exact C*-algebra, then
$\mathcal{A}$ is a von Neumann algebra.

It is straightforward that if $q$ is a quasitrace on some
C*-algebra $\mathfrak{A}$, and $A\in\mathfrak{A}$ is some element that can
be written as $A=XX^*-X^*X$ for some $X\in\mathfrak{A}$,
such that $XX^*$ and $X^*X$ {\em commute}, then $q(A)=0$. In this paper we are
going to take a closer look at such $A$'s, which will be referred to as
{\em abelian self-commutators}.

Suppose now $\mathcal{A}$ is a finite AW*-factor, which is
contained as an AW*-subalgebra in a finite AW*-factor $\mathcal{B}$.
Due to the uniquess of the quasitrace,  for $A\in\mathcal{A}$, one has the equivalence
$q_{\mathcal{A}}(A)=0\Leftrightarrow q_{\mathcal{B}}(A)=0$, so a sufficient
condition for $q_{\mathcal{A}}(A)=0$ is that $A$ is an abelian seff-commutator
in $\mathcal{B}$. In this paper we prove the converse, namely: {\em If
$\mathcal{A}$ is a finite AW*-factor, and $A\in\mathcal{A}$ is a self-adjoint element
of quasitrace zero, then there exists a finite AW*-factor $\mathcal{M}$, that
contains $\mathcal{A}$ as an AW*-subfactor, such that $A$ is an abelian
self-commutator in $\mathcal{M}$}. Moreover, $\mathcal{M}$ can be chosen such that
it is of same type
as $\mathcal{A}$, and {\em every self-adjoint element $X\in\mathcal{M}$ of
quasitrace zero is an abelian self commutator in $\mathcal{M}$}.
Specifically,
in the type \typeone{n}, $\mathcal{M}$ is
$\mathcal{A}$ itself, and in the type \twone\ case, $\mathcal{M}$
is an ultraproduct.

The paper is organized as follows. In Section 1 we introduce our notations,
and we recall several standard
results from the literature, and in Section 2 we prove the main results.

\

\section{Prelimiaries}

\begin{notations}
Let $\mathcal{A}$ be a unital C*-algebra.
\begin{itemize}
\item[A.] We denote by $\mathcal{A}_{sa}$ the real linear space of self-adjoint elements.
We denote by
$\mathbf{U}(\mathcal{A})$
the group of unitaries in $\mathcal{A}$.
We denote by $\mathbf{P}(\mathcal{A})$ the collection of
projections in $\mathcal{A}$, that is,
$\mathbf{P}(\mathcal{A})=\big\{P\in\mathcal{A}_{sa}\,:\,P=P^2\big\}$.
\item[B.] Two elements $A,B\in\mathcal{A}$ are said to be
{\em unitarily equivalent in $\mathcal{A}$}, in which case we write
$A\approx B$, if there exists $U\in\mathbf{U}(\mathcal{A})$ such that
$B=UAU^*$.
\item[C.]
Two elements $A, B\in\mathcal{A}$ are said to be {\em orthogonal},
in which case we write $A\perp B$, if:
$AB=BA=AB^*=B^*A=0$.
(Using the Fuglede-Putnam Theorem, in the case when one of the two is
normal, the above condition reduces to: $AB=BA=0$. If both $A$
and $B$ are normal, one only needs $AB=0$.)
A collection $(A_j)_{j\in J}\subset\mathcal{A}$ is said to be orthogonal, if
$A_i\perp A_j$, $\forall\,i\neq j$.
\end{itemize}
\end{notations}

Finite AW*-factors have several interesting features, contained in the
following well-known result (stated without proof).
\begin{proposition}
Assume $\mathcal{A}$ is a finite AW*-factor.
\begin{itemize}
\item[A.] For any element $X\in\mathcal{A}$, one has: $XX^*\approx X^*X$.
\item[B.] If $X_1,X_2,Y_1,Y_2\in\mathcal{A}$ are such that $X_1\approx X_2$,
$Y_1\approx Y_2$, and $X_k\perp Y_k$, $k=1,2$, then $X_1+Y_1\approx X_2+Y_2$.
\end{itemize}
\end{proposition}

\begin{definition}
Let $\mathcal{A}$ be a unital C*-algebra.
An element $A\in\mathcal{A}_{sa}$ is called an {\em abelian
self-commutator}, if there exists $X\in \mathcal{A}$, such that
\begin{itemize}
\item $(XX^*)(X^*X)=(X^*X)(XX^*)$;
\item $A=XX^*-X^*X$.
\end{itemize}
\end{definition}

\begin{remark}
It is obvious that if $A\in\mathcal{A}_{sa}$ is an abelian
self-commutator, 
then
$q(A)=0$, for any quasitrace $q$ on $\mathcal{A}$.
\end{remark}
\addtocounter{remark}{1}
Abelian self-commutators in finite AW*-factors
can be characterized as follows.

\begin{proposition}
Let $\mathcal{A}$ be a finite AW*-factor. For an element $A\in\mathcal{A}_{sa}$,
the following are equivalent:
\begin{itemize}
\item[(i)] $A$ is an abelian self-commutator in $\mathcal{A}$;
\item[(ii)] there exists $A_1,A_2\in\mathcal{A}_{sa}$ with:
\begin{itemize}
\item[$\bullet$] $A_1A_2=A_2A_1$;
\item[$\bullet$] $A=A_1-A_2$;
\item[$\bullet$] $A_1\approx A_2$.
\end{itemize}
\end{itemize}
\end{proposition}
\begin{proof}
The implication (i)$\Rightarrow$(ii) is trivial by Proposition 1.1.

Conversely, assume $A_1$ and $A_2$ are as in (ii), and let
$U\in\mathbf{U}(\mathcal{A})$ be such that $UA_1U^*=A_2$. Choose a real number
$t>0$, such that $A_1+tI\geq 0$ (for example $t=\|A_1\|$), and define the element
$X=(A_1+tI)^{1/2}U^*$. Notice that $XX^*=A_1+tI$, and $X^*X=A_2+tI$, so
$XX^*$ and $X^*X$ commute. Now we are done, since $XX^*-X^*=A_1-A_2=A$.
\end{proof} 

\begin{notation}
In \cite{Ha} Haagerup shows that, given a normalized quasitrace $q$ on a
unital C*-algebra $\mathcal{A}$, the map $d_q:\mathcal{A}\times
\mathcal{A}\to [0,\infty)$, given by
$$d_q(X,Y)=q\big((X-Y)^*(X-Y))^{\frac 13}, \,\,\,
\forall\,X,Y\in\mathcal{A},$$
defines a metric.
We refer to this metric as the
{\em Haagerup ``$\frac 23$-metric'' associated with $q$}.
Using the inequality $|q(X)|\leq 2\|X\|$, one also has the
inequality
\begin{equation}
d_q(X,Y)\leq\root{3}\of{2}\|X-Y\|^{\frac 23},\,\,\,\forall\,X,Y\in\mathcal{A}.
\label{d-norm}
\end{equation}
If $\mathcal{A}$ is a finite AW*-factor, we denote by $q_{\mathcal{A}}$ the
(unique) normalized quasitrace on $\mathcal{A}$, and we denote by
$d_{\mathcal{A}}$ the Haagerup ``$\frac 23$-metric'' associated with
$q_{\mathcal{A}}$.
\end{notation}

We now concentrate on some issues that
deal with the problem of ``enlarging'' a finite AW*-factor to
a ``nicer'' one.
Recall that, given an AW*-algebra $\mathcal{B}$, a subset
$\mathcal{A}\subset\mathcal{B}$ is declared an {\em AW*-subalgebra of
$\mathcal{B}$}, if it has the following properties:
\begin{itemize}
\item[(i)] $\mathcal{A}$ is a C*-subalgebra of $\mathcal{B}$;
\item[(ii)] $\mathbf{s}(A)\in\mathcal{A}$, $\forall\,A\in\mathcal{A}$;
\item[(iii)] if $(P_i)_{i\in I}\subset\mathbf{P}(\mathcal{A})$, then
$\bigvee_{i\in I}P_i\in\mathcal{A}$.
\end{itemize}
(In condition (ii) the projection $\mathbf{s}(A)$ is the support of $A$ in
$\mathcal{B}$. In (iii) the supremum is computed in $\mathcal{B}$.)
In this case it is pretty clear that $\mathcal{A}$ is an AW*-algebra on its own,
with unit $I_{\mathcal{A}}=\bigvee_{A\in\mathcal{A}_{sa}}\mathbf{s}(A)$.
Below we take a look at the converse statement, namely at the question whether a
C*-subalgebra $\mathcal{A}$ of an AW*-algebra $\mathcal{A}$, which is an AW*-algebra on its
own, is in fact an AW*-subalgebra of $\mathcal{B}$. We are going to restrict
ourselves with the factor case, and for this purpose we introduce the
following terminology.
\begin{definition}
Let $\mathcal{B}$ be an AW*-factor. An AW*-subalgebra
$\mathcal{A}\subset\mathcal{B}$ is called an {\em AW*-subfactor of
$\mathcal{B}$}, if $\mathcal{A}$ is a factor, and $\mathcal{A}\ni I$ -- the unit
in $\mathcal{B}$.
\end{definition}
\begin{proposition}
Let $\mathcal{A}$ and $\mathcal{B}$ be finite AW*-factors. If
$\pi:\mathcal{A}\to\mathcal{B}$ be a unital (i.e. $\pi(I)=I$)
$*$-homomorphism, then $\pi(\mathcal{A})$ is a AW*-subfactor of $\mathcal{B}$.
\end{proposition}
\begin{proof}
Denote for simplicity $\pi(\mathcal{A})$ by $\mathcal{M}$.
Since $\mathcal{A}$ is simple, $\pi$ is injective, so $\mathcal{M}$ is
$*$-isomprhic to $\mathcal{A}$. Among other things, this shows that
$\mathcal{M}$ is a factor, which contains the unit $I$ of $\mathcal{B}$.
We now proceed to check the two key
conditions (ii) and (iii) that ensure that $\mathcal{M}$ is an AW*-subalgebra in $\mathcal{B}$.

(ii). Start with some element $M\in\mathcal{M}$, written as $M=\pi(A)$, for some
$A\in\mathcal{A}$, and let us show that $\mathbf{s}(M)$ -- the support of $M$ in
$\mathcal{B}$ --
in fact belongs to $\mathcal{M}$. This will be the result of the
following:
\begin{claimone}
One has the equality 
$\mathbf{s}(M)=\pi\big(\mathbf{s}(A)\big)$, where
$\mathbf{s}(A)$ denotes the support of $A$ in $\mathcal{A}$.
\end{claimone}
Denote the projection $\pi\big(\mathbf{s}(A)\big)\in\mathbf{P}(\mathcal{M})$
by $P$. First of all, since $\big(I-\mathbf{s}(A)\big)A=0$ (in $\mathcal{A}$),
we have $(I-P)M=0$ in $\mathcal{B}$, so $(I-P)\perp\mathbf{s}(M)$, i.e.
$\mathbf{s}(M)\geq P$. Secondly, since $q_{\mathcal{B}}\circ\pi:\mathcal{A}\to
\mathbb{C}$ is a quasitrace, we must have the equality
\begin{equation}
q_{\mathcal{B}}\circ\pi=q_{\mathcal{A}}.
\label{qcomp}
\end{equation}
In particular the projection $P$ has dimension
$D_{\mathcal{B}}(P)=D_{\mathcal{A}}\big(\mathbf{s}(A)\big)$.
We know however that for a self-adjoint element
$X$ in a finite AW*-factor with quasitrace $q$, one has the equality
$q\big(\mathbf{s}(X)\big)=\mu^X\big(\mathbb{R}\smallsetminus\{0\}\big)$, where
$\mu^X$ is the scalar spectral measure, defined implictly (using Riesz' Theorem) by
$$\int_{\mathbb{R}}f\,d\mu^X=q\big(f(X)\big),\,\,\,\forall\,f\in C_0(\mathbb{R}).$$
So in our case we have the equalities
\begin{align}
D_{\mathcal{B}}\big(\mathbf{s}(M)\big)&=\mu^M_{\mathcal{B}}\big(\mathbb{R}\smallsetminus\{0\}\big),
\label{DsM}\\
D_{\mathcal{B}}\big(P)&=D_{\mathcal{A}}\big(\mathbf{s}(A)\big)=\mu^A_{\mathcal{A}}\big(\mathbb{R}\smallsetminus\{0\}\big),
\label{DP}
\end{align}
where the subscripts indicate the ambient AW*-factor.
Since $\pi$ is a $*$-homoorphism, one has the equality $\pi\big(f(A)\big)=f(M)$,
$\forall\,f\in C_0(\mathbb{R})$, and then by \eqref{qcomp} we get
$$
\int_{\mathbb{R}}f\,d\mu^M_{\mathcal{B}}=q_{\mathcal{B}}\big(f(M)\big)=
(q_{\mathcal{B}}\circ\pi)\big(f(A)\big)=
q_{\mathcal{A}}\big(f(A)\big)=
\int_{\mathbb{R}}f\,d\mu^A_{\mathcal{A}},
\,\,\,\forall\,f\in C_0(\mathbb{R}).$$
In particular we have the equality $\mu^M_{\mathcal{B}}=\mu^A_{\mathcal{A}}$,
and then \eqref{DsM} and \eqref{DP} will force $D_{\mathcal{B}}\big(\mathbf{s}(M)\big)=
D_{\mathcal{B}}(P)$. Since $P\geq\mathbf{s}(M)$, the equality of dimensions will force
$P=\mathbf{s}(M)$.

(iii). Start with a collection of projections $(P_i)_{i\in I}\subset\mathbf{P}(\mathcal{M})$, 
let $P=\bigvee_{i\in I}P_i$ (in $\mathcal{B}$), and let us prove that $P\in\mathcal{M}$.
Write each $P_i=\pi(Q_i)$, with $Q_i\in\mathbf{P}(\mathcal{A})$, and let $Q=\bigvee_{i\in I}Q_i$
(in $\mathcal{A}$). The desired conclusion will result from the following.
\begin{claimtwo}
$P=\pi(Q)$.
\end{claimtwo}
Denote by $\mathcal{F}$ the collection of all finite subsets of $I$, which
becomes a directed set with inclusion, and define the nets
$P_F=\bigvee_{i\in F}P_i$ (in $\mathcal{B}$) and $Q_F=\bigvee_{i\in F}Q_i$ (in $\mathcal{A}$).
On the one hand, if we consider the element $X_F=\sum_{i\in F}Q_i$, then
$Q_F=\mathbf{s}(X_F)$ (in $\mathcal{A}$), and $P_F=\mathbf{s}\big(\sum_{i\in I}P_i\big)=
\mathbf{s}\big(\pi(X_F)\big)$ (in $\mathcal{B}$), so by Claim 1, we have the equality
$P_F=\pi (Q_F)$. On the other hand, we have $Q=\bigvee_{F\in\mathcal{F}}Q_F$ (in $\mathcal{A}$),
with the net $(Q_F)_{F\in\mathcal{F}}$ increasing, so we get the equality
$D_{\mathcal{A}}(Q)=\lim_{F\in\mathcal{F}}D_{\mathcal{A}}(Q_F)$.
Arguing the same way (in $\mathcal{B}$), and using the equalities
$P_F=\pi(Q_F)$, we get
$$
D_{\mathcal{B}}(P)=\lim_{F\in\mathcal{F}}D_{\mathcal{B}}(P_F)=
\lim_{F\in\mathcal{F}}D_{\mathcal{B}}\big(\pi(Q_F)\big)=
\lim_{F\in\mathcal{F}}D_{\mathcal{A}}(Q_F)=
D_{\mathcal{A}}(Q)
=D_{\mathcal{B}}\big(\pi(Q)\big).$$
Finally, since $\big[I-\pi(Q)\big]P_i=\pi\big([I-Q]Q_i\big)=0$, $\forall\,i\in
I$, we get the inequality $\pi(Q)\geq P$, and then the equality
$D_{\mathcal{B}}(P)=D_{\mathcal{B}}\big(\pi(Q)\big)$ will force $P=\pi(Q)$.
\end{proof}

\begin{comment}
In the above proof we employed an argument based on the following property of the
dimension function $D$ on a finite AW*-factor $\mathcal{A}$:
\begin{itemize}
\item[{\sc (l)}]
{\em If a net
$(P_\lambda)_{\lambda\in\Lambda}\subset\mathbf{P}(\mathcal{A})$
is increasing, then $D\big(\bigvee_{\lambda\in\Lambda}P_\lambda\big)=\lim_{\lambda\in\Lambda}D(P_\lambda)$}.
\end{itemize}
The literature (\cite{Be},\cite{Kap}) often mentions a different feature: the
{\em complete additivity}
\begin{itemize}
\item[{\sc (c.a.)}]
{\em If $(E_i)_{i\in I}\subset\mathbf{P}(\mathcal{A})$ is an orthogonal family,
then $D\big(\bigvee_{i\in I}E_i\big)=\sum_{i\in I}D(E_i)$}.
\end{itemize}
To prove property {\sc (l)} one can argue as follows. Let
$P=\bigvee_{\lambda\in\Lambda}P_\lambda$, so that $D(P)\geq D(P_\lambda)$,
$\forall\,\lambda\in\Lambda$. In particular if we take
$\ell=\lim_{\lambda\in\Lambda}D(P_\lambda)$ (which exists by monotonicity), we get
$D(P)\geq\ell$. To prove that in fact we have $D(P)=\ell$, we construct a sequence
$\alpha_1\prec\alpha_2\prec\dots$ in $\Lambda$, such that $\lim_{n\to\infty}
D(P_{\alpha_n})=\ell$, we consider the projection
$Q=\bigvee_{n\in\mathbb{N}}P_{\alpha_n}\leq P$, and we show first that
$D(Q)=\ell$, then $Q=P$. The equality $D(Q)=\ell$
follows from {\sc (c.a.)} since now
we can write
$Q=P_{\alpha_1}\vee\bigvee_{n=1}^\infty [P_{\alpha_{n+1}}-P_{\alpha_n}]$ with
all projections orthogonal, so we get
\begin{gather*}
D(Q)=D(P_{\alpha_1})+\sum_{n=1}^\infty D(P_{\alpha_{n+1}}-P_{\alpha_n})=\\
=D(P_{\alpha_1})+\sum_{n=1}^\infty \big[D(P_{\alpha_{n+1}})-D(P_{\alpha_n})\big]=
\lim_{n\to\infty}D(P_{\alpha_n})=\ell.\end{gather*}
To prove the equality $Q=P$, we fix for the moment $\lambda\in\Lambda$,
and integer $n\geq 1$, and some
$\mu\in\Lambda$ with $\mu\succ\lambda$ and $\mu\succ\alpha_n$, and we observe
that, using the Paralelogram Law, we have:
$$
P_\lambda-P_\lambda\wedge Q\leq
P_\lambda-P_\lambda\wedge P_{\alpha_n}\approx
P_\lambda\vee P_{\alpha_n}-P_{\alpha_n}\leq P_\mu-P_{\alpha_n},$$
so applying the dimension function we get
$$D(P_\lambda-P_\lambda\wedge Q)\leq D(P_\mu-P_{\alpha_n})=
D(P_\mu)-D(P_{\alpha_n})\leq\ell-D(P_{\alpha_n}).$$
Since the inequality
$D(P_\lambda-P_\lambda\wedge Q)\leq\ell-D(P_{\alpha_n})$ holds for arbitrary
$n\in\mathbb{N}$ and $\lambda\in\Lambda$, taking limit (as $n\to\infty$) yields
$D(P_\lambda-P_\lambda\wedge Q)=0$, which in turn forces
$P_\lambda=P_\lambda\wedge Q$, which means that $P_\lambda\leq Q$. Since this is
true for all $\lambda\in\Lambda$, it will force $Q\geq P$, so we must have $Q=P$.
\end{comment}

We now recall the ultraproduct construction of finite AW*-factors, discussed
for example in
\cite{BH} and \cite{Ha}.
\begin{notations}
Let $\mathbf{A}=(\mathcal{A}_n)_{n\in\mathbb{N}}$ be a sequence of finite
AW*-factors, and let $q_n:\mathcal{A}_n\to\mathbb{C}$ denote the (unique)
normalized quasitrace on $\mathcal{A}_n$. One considers the finite
AW*-algebra
$$\mathbf{A}^\infty=\big\{(X_n)_{n\in\mathbb{N}}\in\prod_{n\in\mathbb{N}}
\mathcal{A}_n\,:\,
\sup_{n\in\mathbb{N}}\|X_n\|<\infty\big\}.$$
Given a free ultrafilter $\mathcal{U}$ on $\mathbb{N}$, one defines the
quasitrace
$\tau_{\mathcal{U}}:\mathbf{A}^\infty\to\mathbb{C}$ by
$$\tau_{\mathcal{U}}(\boldsymbol{x})=\lim_{\mathcal{U}}q_n(X_n),\,\,\,\forall\,
\boldsymbol{x}=(X_n)_{n\in\mathbb{N}}\in\mathbf{A}^\infty.$$
Next one considers the norm-closed ideal
$$\mathbf{J}_{\mathcal{U}}=\big\{\boldsymbol{x}\in
\mathbf{A}^\infty\,:\,\tau_{\mathcal{U}}(\boldsymbol{x}^*\boldsymbol{x})=0
\big\}.$$
It turns out that quotient C*-algebra
$\mathbf{A}_{\mathcal{U}}=\mathbf{A}^\infty/\mathbf{J}_{\mathcal{U}}$ becomes a
finite AW*-factor. Moreover, its (unique) normalized quasitrace
$q_{\mathbf{A}_{\mathcal{U}}}$ is defined implictly by
$q_{\mathbf{A}_{\mathcal{U}}}\circ\Pi_{\mathcal{U}}=\tau_{\mathcal{U}}$, where
$\Pi_{\mathcal{U}}:\mathbf{A}^\infty\to\mathbf{A}_{\mathcal{U}}$ denotes the
quotient $*$-homomorphism.
\end{notations}

The finite AW*-factor $\mathbf{A}_{\mathcal{U}}$ is referred to as the
{\em $\mathcal{U}$-ultraproduct of the sequence $\mathbf{A}=
(\mathcal{A}_n)_{n\in\mathbb{N}}$}.
\begin{remarks} Let
$\mathbf{A}=(\mathcal{A}_n)_{n\in\mathbb{N}}$ be a sequence
of finite factors.
\begin{itemize}
\item[A.]
With the notations above, if $\boldsymbol{x}=(X_n)_{n\in\mathbb{N}},\boldsymbol{y}=(Y_n)_{n\in\mathbb{N}}\in\mathbf{A}^\infty$
are elements that satisfy the condition
$\lim_{\mathcal{U}}d_{\mathcal{A}_n}(X_n,Y_n)=0$,
then $\Pi_{\mathcal{U}}(\boldsymbol{x})=\Pi_{\mathcal{U}}(\boldsymbol{y})$ in
$\mathbf{A}_{\mathcal{U}}$.
This is trivial, since the given condition forces
$$\lim_{\mathcal{U}}q_{\mathcal{A}_n}\big((X_n-Y_n)^*(X_n-Y_n)\big)=0,$$
i.e.
$\boldsymbol{x}-\boldsymbol{y}\in\mathbf{J}_{\mathcal{U}}$.
\item[B.] For $\boldsymbol{x}=(X_n)_{n\in\mathbb{N}}\in\mathbf{A}^\infty$, one
has the inequality:
$\|\Pi_{\mathcal{U}}(\boldsymbol{x})\|\leq\lim_{\mathcal{U}}\|X_n\|$.
To prove this inequality we start off by denoting
$\lim_{\mathcal{U}}\|X_n\|$ by $\ell$, and we observe that
given any $\varepsilon>0$, the set
$$U_\varepsilon=\{n\in\mathbb{N}\,:\,\ell-\varepsilon<\|X_n\|<\ell+\varepsilon\}$$
belongs to $\mathcal{U}$, so if we define
the sequence
$\boldsymbol{x}_\varepsilon=(X^\varepsilon_n)_{n\in\mathbb{N}}$ by
$$X^\varepsilon_n=\left\{\begin{array}{cl} X_n&\text{ if }
n\not\in U_\varepsilon\\
0&\text{ if }n\in U_\varepsilon
\end{array}\right.$$
we clearly have $\lim_{\mathcal{U}}\|X_n^\varepsilon\|=0$.
In particular, by part A, we have
$\Pi_{\mathcal{U}}(\boldsymbol{x})=\Pi_{\mathcal{U}}(\boldsymbol{x}-
\boldsymbol{x}_\varepsilon)$. Since $\|X_n-X_n^\varepsilon\|\leq
\ell+\varepsilon$, $\forall\,n\in\mathbb{N}$, it follows that
$$\|\Pi_{\mathcal{U}}(\boldsymbol{x})\|=
\|\Pi_{\mathcal{U}}(\boldsymbol{x}-
\boldsymbol{x}_\varepsilon)\|\leq\ell+\varepsilon,$$
and since the inequality $\|\Pi_{\mathcal{U}}(\boldsymbol{x})\|\leq
\ell+\varepsilon$ holds for all $\varepsilon>0$, it follows that
we indeed have
$\|\Pi_{\mathcal{U}}(\boldsymbol{x})\|\leq\ell$.
\end{itemize}
\end{remarks}
\begin{example}
Start with a finite AW*-factor $\mathcal{A}$ of type \twone\ and a free ultrafilter
$\mathcal{U}$ on $\mathbb{N}$. Let $\mathcal{A}_{\mathcal{U}}$ denote the
ultrapoduct of the constant sequence $\mathcal{A}_n=\mathcal{A}$.
For every $X\in\mathcal{A}$ let
$\Gamma(X)=(X_n)_{n\in\mathbb{N}}\in\mathbf{A}^\infty$ be the constant sequence:
$X_n=X$. It is obvious that $\Gamma:\mathcal{A}\to\mathbf{A}^\infty$ is a
unital $*$-homomorphism, so the composition $\Delta_{\mathcal{U}}=
\Pi_{\mathcal{U}}\circ\Gamma
:\mathcal{A}\to\mathcal{A}_{\mathcal{U}}$ is again a
unital
$*$-homomorphism. Using Proposition 2.1 it follows that
$\Delta_{\mathcal{U}}(\mathcal{A})$ is an AW*-subfactor in
$\mathcal{A}_{\mathcal{U}}$.

Moreover, if $\mathcal{B}$ is some finite AW*-factor, and
$\pi:\mathcal{B}\to\mathcal{A}$ is some unital
$*$-homomorphism, then the $*$-homomorphism
$\pi_{\mathcal{U}}=\Delta_{\mathcal{U}}\circ\pi:\mathcal{B}\to\mathcal{A}_{\mathcal{U}}$
gives rise to an AW*-subfactor $\pi_{\mathcal{U}}(\mathcal{B})$ of
$\mathcal{A}_{\mathcal{U}}$.
\end{example}

\

\section{Main Results}

We start off with the analysis of the type \onefin\ situation, i.e.
the algebras of the form $\text{\rm Mat}_n(\mathbb{C})$ -- the $n\times n$ complex
matrices.
To make the exposition a little easier, we are going to use the un-normalized trace
$\tau:\text{\rm Mat}_n(\mathbb{C})\to \mathbb{C}$ with $\tau (I_n)=n$.

The main result in the type \onefin\ --
stated in a way that will allow an inductive proof --
is as follows.
\begin{theorem}
Let $n\geq 1$ be an integer, and let $X\in \text{\rm Mat}_n(\mathbb{C})_{sa}$ be a
matrix with $\tau(X)=0$.
\begin{itemize}
\item[A.] For any projection $P\in \text{\rm Mat}_n(\mathbb{C})$ with
$PX=XP$ and $\tau(P)=1$, there exist elements
$A,B\in \text{\rm Mat}_n(\mathbb{C})_{sa}$ with:
\begin{itemize}
\item[$\bullet$] $AB=BA$;
\item[$\bullet$] $X=A-B$;
\item[$\bullet$] $A\approx B$;
\item[$\bullet$] $\max\big\{\|A\|,\|B\|\big\}\leq \|X\|$;
\item[$\bullet$] $A\perp P$.
\end{itemize}
\item[B.] $X$ is an abelian self-commutator in $\text{\rm Mat}_n(\mathcal{C})$.
\end{itemize}
\end{theorem}
\begin{proof}
A. We are going to use induction on $n$.
The case $n=1$ is trivial, since it forces $X=0$, so we can take
$A=B=0$.
Assume now property A is true for all $n<N$, and let us prove it for
$n=N$. Fix some $X\in \text{\rm Mat}_N(\mathbb{C})$ with $\tau(X)=0$, and
a projection $P\in \text{\rm Mat}_N(\mathbb{C})$ with $\tau(P)=1$ such that
$PX=XP$. The case $X=0$ is trivial, so we are going to assume $X\neq 0$.
Let us list the spectrum of $X$ as
$\text{Spec}(X)=\{\alpha_1<\alpha_2<\dots<\alpha_m\}$, and let $(E_i)_{i=1}^m$
be the corresponding spectral projections, so that
\begin{itemize}
\item[(i)] $\tau(E_i)>0$, $\forall\,i\in\{1,\dots,m\}$;
\item[(ii)] $E_i\perp E_j$, $\forall\,i\neq j$, and $\sum_{i=1}^mE_i=I_N$;
\item[(iii)] $X=\sum_{i=1}^m\alpha_i E_i$, so
$\tau(X)=\sum_{i=1}^m\alpha_i\tau(E_i)$.
\end{itemize}
Since $\tau(P)=1$ and $P$ commutes with $X$, there exists a unique index $i_0\in\{1,\dots,m\}$
such that $P\leq E_{i_0}$. Since none of the
inclusions $\text{Spec}(X)\subset (0,\infty)$ or
$\text{Spec}(X)\subset (-\infty,0)$ is possible, there exists
$i_1\in\{1,\dots,m\}$, $i_1\neq i_0$,
such that one of the following inequalities holds
\begin{gather}
\alpha_{i_1}<0\leq\alpha_{i_0},\label{1le0}\\
\alpha_{i_1}>0\geq\alpha_{i_0}.\label{1ge0}\end{gather}
Choose then a projection $Q\leq E_{i_1}$ with $\tau(Q)=1$, and let us define the
elements $S=\alpha_{i_0}(P-Q)$ and $Y=X-S$.
Notice that
$$Y=\sum_{i\neq i_0,i_0}\alpha_iE_i+\alpha_{i_1}(E_{i_1}-Q)+\alpha_{i_0}(E_{i_0}-P)+
(\alpha_{i_1}+\alpha_{i_0})Q,$$
so in particular we have $Y\perp P$. Notice also that
either one of \eqref{1le0} or \eqref{1ge0} yields
$$|\alpha_{i_1}+\alpha_{i_0}|\leq\max\big\{|\alpha_{i_0}|,|\alpha_{i_1}|\big\}\leq\|X\|,$$
so we have $\|Y\|\leq\|X\|$. Finally, since both
$Y$ an $Q$ belong to the subalgebra
$$\mathcal{A}=(I_N-P)\text{\rm Mat}_N(\mathbb{C})(I_N-P),$$
which is
$*$-isomorphic to $\text{\rm Mat}_{N-1}(\mathbb{C})$, using the inductive hypothesis, with
$Y$ and $Q$ (which obviously commute), there exist
$A_0,B_0\in \mathcal{A}$, with
\begin{itemize}
\item $A_0B_0=B_0A_0$;
\item $Y=A_0-B_0$;
\item $A_0\approx B_0$;
\item $\max\big\{\|A_0\|,\|B_0\|\big\}\leq \|Y\|\leq\|X\|$;
\item $A_0\perp Q$.
\end{itemize}
It is now obvious that the elements $A=A_0-\alpha_{i_0}Q$ and
$B=B-\alpha_{i_0}P$ will satisfy the desired hypothesis (at one point,
Proposition 1.1.B is invoked).

B. This statement is obvious from part A, since one can always start with an
arbitrary projection $P\leq E_m$, with $\tau(P)=1$, and such a projection
obviously commutes with $X$.
\end{proof}

In preparation for the type \twone\ case, we have the following approximation result.
\begin{lemma}
Let $\mathcal{A}$ be an AW*-factor of type \twone, and let
$\varepsilon>0$ be a real number. For any element $X\in\mathcal{A}_{sa}$, there
exists an AW*-subfactor $\mathcal{B}\subset\mathcal{A}$, of type
\onefin, and an element $B\in\mathcal{B}_{sa}$ with
\begin{itemize}
\item[(i)] $d_{\mathcal{A}}(X,B)<\varepsilon$;
\item[(ii)] $q_{\mathcal{A}}(X)=q_{\mathcal{A}}(B)$;
\item[(iii)] $\|B\|\leq\|X\|+\varepsilon$.
\end{itemize}
\end{lemma}
\begin{proof}
We begin with the following
\begin{partcase}
Assume $X$ has finite spectrum.
\end{partcase}
Let $\text{Spec}(X)=\{\alpha_1<\alpha_2<\dots<\alpha_m\}$, and let
$E_1,\dots,E_m$ be the corresponding spectral projections, so that
\begin{itemize}
\item $D(E_i)>0$, $\forall\,i\in\{1,\dots,m\}$;
\item $E_i\perp E_j$, $\forall\,i\neq j$, and $\sum_{i=1}^mE_i=I$;
\item $X=\sum_{i=1}^m\alpha_i E_i$, so
$q_{\mathcal{A}}(X)=\sum_{i=1}^m\alpha_i D(E_i)$.
\end{itemize}

For any integer $n\geq 2$ define the
set $Z_n=\big\{0,\frac 1n,\frac 2n,\dots,\frac{n-1}n,1\big\}$, and let
$\theta_n:\{1,\dots,m\}\to Z_n$ be the map defined by
$$\theta_n(i)=\max\big\{\zeta\in Z_n\,:\,\zeta\leq D(E_i)\big\}.$$

For every $i\in\{1,\dots,m\}$, and every integer $n\geq 2$, chose
$P_{ni}\in\mathbf{P}(\mathcal{A})$ be an arbitrary
projection with $P_{ni}\leq E_i$, and
$D(P_{ni})=\theta_n(i)$.
Notice that, for a fixed $n\geq 2$, the projections
$P_{n1},\dots,P_{nm}$ are pairwise orthogonal, and have dimensions
in the set $Z_n$, hence there exists a subfactor $\mathcal{B}_n$ of type
\typeone{n}, that contains them. Define then the element $H_n=\sum_{i=1}^m
\alpha_iP_{ni}\in(\mathcal{B}_n)_{sa}$.
Note that
$\|H_n\|\leq \|X\|$.
We wish to prove that
\begin{itemize}
\item[{\sc (a)}] $\lim_{n\to\infty}d_{\mathcal{A}}(X,H_n)=0$;
\item[{\sc (b)}] $\lim_{n\to\infty}q_{\mathcal{A}}(H_n)=q_{\mathcal{A}}(X)$.
\end{itemize}
To prove these assertions, we first observe that, for each $n\geq 2$, the
elements $X$ and $H_n$ commute, and we have
$$X-H_n=\sum_{i=1}^m\alpha_i(E_i-P_{ni}).$$
In particular, one has
\begin{equation}
|q_{\mathcal{A}}(X)-q_{\mathcal{A}}(H_n)|\leq
\sum_{i=1}^m|\alpha_i|\cdot D(E_i-P_{ni})\leq
m\|X\|\cdot\max\big\{D(E_i-P_{ni})\}_{i=1}^m.
\label{qXH}
\end{equation}
Likewise, since
$$(X-H_n)^*(X-H_n)=
\sum_{i=1}^m\alpha_i^2(E_i-P_{ni}),$$
we have
\begin{equation}
\begin{split}
&q_{\mathcal{A}}\big((X-H_n)^*(X-H_n)\big)=
\sum_{i=1}^m\alpha_i^2\cdot D(E_i-P_{ni})\leq\\
&\qquad\qquad\leq m\|X\|^2\cdot\max\big\{D(E_i-P_{ni})\big\}_{i=1}^m.
\end{split}
\label{qXH2}
\end{equation}
By construction however we have $D(E_i-P_{ni})<\frac 1n$, so
the estimates \eqref{qXH} and \eqref{qXH2} give
\begin{gather*}
|q_{\mathcal{A}}(X)-q_{\mathcal{A}}(H_n)|\leq\frac{m\|X\|}n\\
d_{\mathcal{A}}(X,H_n)\leq\root{3}\of{\frac{m\|X\|^2}n},
\end{gather*}
which clearly give the desired assertions {\sc (a)} and {\sc (b)}.

Using the conditions {\sc (a)} and {\sc (b)}, we immediately see that, if we
define the numbers $\beta=q_{\mathcal{A}}(X)$,
$\beta_n=q_{\mathcal{A}}(H_n)$, and the elements
$B_n=H_n+(\beta-\beta_n)I\in\mathcal{B}_n$,
then the sequence $(B_n)_{n\geq 2}$ will still satisfy
$\lim_{n\to\infty}d_{\mathcal{A}}(X,B_n)=0$, but also $q_{\mathcal{A}}(B_n)=
q_{\mathcal{A}}(B_n)$, and
$$\|B_n\|\leq\|H_n\|+|\beta-\beta_n\|\leq\|X\|+|\beta-\beta_n|,$$
which concludes the proof of the Particular Case.

Having proven the Particular Case, we now proceed with the general case.
Start with an arbitrary element $X\in\mathcal{A}_{sa}$, and pick a
sequence
$(T_n)_{n\in\mathbb{N}}\subset\mathcal{A}_{sa}$ of elements with finite
spectrum, such that $\lim_{n\to\infty}\|T_n-X\|=0$.
(This can be done using Borel functional calculus.)
Using the norm-continuity of the quasitrace, we have
$\lim_{n\to\infty}q_{\mathcal{A}}(T_n)=q_{\mathcal{A}}(X)$, so if we define
$X_n=T_n+\big(q_{\mathcal{A}}(X)-q_{\mathcal{A}}(T_n)\big)I$, we will still have
$\lim_{n\to\infty}\|X_n-X\|=0$, but also $q_{\mathcal{A}}(X_n)=q_{\mathcal{A}}(X)$.
In particular, there exists some $k\geq 1$, such that
\begin{itemize}
\item $d_{\mathcal{A}}(X_k,X)<\varepsilon/2$;
\item $\|X_k\|<\|X\|+\varepsilon/2$.
\end{itemize}
Finally, applying the Particular Case, we can also find an AW*-subfactor
$\mathcal{B}\subset\mathcal{A}$, of type \onefin, and an element
$B\in\mathcal{B}_{sa}$ with
\begin{itemize}
\item $d_{\mathcal{A}}(X_k,B)<\varepsilon/2$;
\item $\|B\|<\|X_k\|+\varepsilon/2$;
\item $q_{\mathcal{A}}(X_k)=q_{\mathcal{A}}(B)$.
\end{itemize}
It is then trivial that $B$ satisfies conditions (i)-(iii).
\end{proof}

We are now in position to prove the main result in the type
\twone\ case.
\begin{theorem}
Let $\mathbf{A}=(\mathcal{A}_n)_{n\in\mathbb{N}}$ be a sequence
of AW*-factors of type \twone,  let $\mathcal{U}$ be a free
ultrafilter on $\mathbb{N}$, and let
$X\in(\mathbf{A}_{\mathcal{U}})_{sa}$ be an element with
$q_{\mathbf{A}_{\mathcal{U}}}(X)=0$. Then
there exist elements
$A,B\in(\mathbf{A}_{\mathcal{U}})_{sa}$ with
\begin{itemize}
\item $AB=BA$;
\item $X=A-B$;
\item $A\approx B$;
\item $\max\big\{\|A\|,\|B\|\big\}\leq\|X\|$.
\end{itemize}
In particular, $X$ is an abelian self-commutator in
$\mathbf{A}_{\mathcal{U}}$.
\end{theorem}
\begin{proof}
Write
$X=\Pi_{\mathcal{U}}(\boldsymbol{x})$, where
$\boldsymbol{x}=(X_n)_{n\in\mathbb{N}}\in\mathbf{A}^\infty$.
Without loss of generality, we can assume that
 all $X_n$'s are self-adjoint, and have norm $\leq \|X\|$.

Consider the elements $\tilde{X}_n=X_n-q_{\mathcal{A}_n}(X_n)I\in
\mathcal{A}_n$. Remark that, since $X_n$ is self-adjoint, we have
$|q_{\mathcal{A}_n}(X_n)\leq
\|X_n\|\leq\|X\|$, so we have $\|\tilde{X}_n\|\leq 2\|X\|$,
$\forall\,n\in\mathbb{N}$,
hence the sequence
$\tilde{\boldsymbol{x}}=(\tilde{X}_n)_{n\in\mathbb{N}}$ defines
an element in
$\mathbf{A}^\infty$. By construction, we have
$\lim_{\mathcal{U}}q_{\mathcal{A}_n}(X_n)=0$, and
$d_{\mathcal{A}_n}(\tilde{X}_n-X_n)=\big|q_{\mathcal{A}_n}(X_n)\big|^{\frac 23}$,
so by Remark 1.2.A it follows that
$X=\Pi_{\mathcal{U}}(\boldsymbol{x})=\Pi_{\mathcal{U}}
(\tilde{\boldsymbol{x}})$.
Note also that
\begin{equation}
\|\tilde{X}_n\|\leq \|X\|+|q_{\mathcal{A}_n}(X_n)|,\,\,\,
\forall\,n\in\mathbb{N}.
\label{norm-tilde}
\end{equation}

Use Lemma 2.1 to find, for each $n\in\mathbb{N}$,
an AW*-subfactor $\mathcal{B}_n$ of $\mathcal{A}_n$, and
and elements
$Y_n\in(\mathcal{B}_n)_{sa}$ with
\begin{itemize}
\item[(i)] $d_{\mathcal{A}_n}(Y_n,\tilde{X}_n)<\frac 1n$;
\item[(ii)] $q_{\mathcal{A}_n}(Y_n)=0$, $\forall\,n\in\mathbb{N}$;
\item[(iii)] $\|Y_n\|\leq\|\tilde{X}_n\|+\frac 1n$.
\end{itemize}
Furthermore, using Theorem 2.1, for each $n\in\mathbb{N}$,
combined with the fact that $q_{\mathcal{B}_n}=q_{\mathcal{A}_n}
\big|_{\mathcal{B}_n}$
(which implies the equality $q_{\mathcal{B}_n}(Y_n)=0$),
there exist
elements $A_n,B_n\in(\mathcal{B}_n)_{sa}$ such that
\begin{itemize}
\item[{\sc (a)}] $A_nB_n=B_nA_n$;
\item[{\sc (b)}] $Y_n=A_n-B_n$;
\item[{\sc (c)}] $A_n\approx B_n$;
\item[{\sc (d)}] $\max\big\{\|A_n\|,\|B_n\|\big\}\leq\|Y_n\|$.
\end{itemize}
Choose $U_n\in\mathbf{U}(\mathcal{B}_n)$, such that $B_n=U_nA_nU_n^*$.

Let us view the sequences
$\boldsymbol{a}=(A_n)_{n\in\mathbb{N}}$, $\boldsymbol{b}=(B_n)_{n\in\mathbb{N}}$,
$\boldsymbol{u}=
(U_n)_{n\in\mathbb{N}}$ as elements in the
AW*-algebra $\mathbf{A}^\infty$, and let us define
the
elements $A=\Pi_{\mathcal{U}}(\boldsymbol{a})$,
$B=\Pi_{\mathcal{U}}(\boldsymbol{b})$, and
$U=\Pi_{\mathcal{U}}(\boldsymbol{u})$ in $\mathcal{A}_{\mathcal{U}}$.
Obviously $A$ and $B$ are self-adjoint. Since $\boldsymbol{u}$ is unitary
in $\mathbf{A}^\infty$, it follows that $U$ is unitary in $\mathcal{A}_{\mathcal{U}}$.
Moreover, since by construction we have
$\boldsymbol{u}\boldsymbol{a}\boldsymbol{u}^*=\boldsymbol{b}$, we also have the equality
$UAU^*=B$, so
$A\approx B$ in $\mathcal{A}$. Finally, since by construction we also
have $\boldsymbol{a}\boldsymbol{b}=\boldsymbol{b}\boldsymbol{a}$, we also get the
equality $AB=BA$.
Since by condition {\sc (d)} we have
$$\max\big\{\|A_n\|,\|B_n\|\big\}
\leq\|Y_n\|\leq\|\tilde{X}_n\|+\tfrac 1n\leq
\|X\|+|q_{\mathcal{A}_n}(X_n)|+\tfrac 1n,$$
by Remark 1.2.B (combined with $\lim_{\mathcal{U}}q_{\mathcal{A}_n}(X_n)
=0$), we get
the
inequality
$$\max\big\{\|A\|,\|B\|\big\}\leq\|X\|.$$
The proof of the Theorem will then be finished, once we prove the
equality $X=A-B$. For this purpose, we consider the
sequences $\tilde{\boldsymbol{x}}=(\tilde{X}_n)_{n\in\mathbb{N}}$ and
$\boldsymbol{y}=(Y_n)_{n\in\mathbb{N}}$, both viewed as elements in
$\mathbf{A}^\infty$. On the one hand, since by construction we
have $\boldsymbol{y}=\boldsymbol{a}-
\boldsymbol{b}$, we get the equality $\Pi_{\mathcal{U}}(\boldsymbol{y})=A-B$.
On the other hand, since
$\lim_{n\to\infty}d_{\mathcal{A}}(Y_n,\tilde{X}_n)=0$, we also
have
$\lim_{\mathcal{U}}d_{\mathcal{A}}(Y_n,\tilde{X}_n)=0$, so
by
Remark 1.2.A we get the equalities
$X=\Pi_{\mathcal{U}}(\boldsymbol{x})=
\Pi_{\mathcal{U}}(\tilde{\boldsymbol{x}})=
\Pi_{\mathcal{U}}(\boldsymbol{y})$, i.e.
$X=A-B$.
\end{proof}

\begin{comment}
Assume $\mathcal{A}$ an AW*-factor of type \twone, and let
$\mathcal{U}$ be a free ultrafilter on $\mathbb{N}$.
Following Example 1.1, $\mathcal{A}$ is identified with the
AW*-subfactor $\Delta_{\mathcal{U}}(\mathcal{A})$ of
$\mathcal{A}_{\mathcal{U}}$.
Under this identification, by Theorem 2.2, every element
$A\in\mathcal{A}_{sa}$ of quasitrace zero is an abelian
self-commutator in $\mathcal{A}_{\mathcal{U}}$.

In connection with this observation, it is legimitate to ask whether $A$
is
in fact an abelian self-commutator in $\mathcal{A}$ itself. The discussion
below aims at answering this question in a somewhat different spirit, based
on the
results from \cite{KN}.
\end{comment}

\begin{definition}
Let $\mathcal{A}$ be an AW*-factor of type \twone. An element
$A\in\mathcal{A}_{sa}$ is called an {\em abelian approximate self-commutator
in $\mathcal{A}$},
if there exist commuting elements $A_1,A_2\in\mathcal{A}_{sa}$ with
$A=A_1-A_2$, and such that $A_1$ and $A_2$ are approximately unitary equivalent,
i.e. there exists a squence $(U_n)_{n\in\mathbb{N}}$ of unitaries in
$\mathcal{A}$ such that $\lim_{n\to\infty}\|U_nA_1U_n^*-A_2\|=0$.
By \cite[Theorem 2.1]{KN} the condition that $A_1$ and $A_2$ are approximately
unitary equivalent -- denoted by $A_1\sim A_2$ -- is equivalent to the condition
$q_{\mathcal{A}}(A_1^k)=q_{\mathcal{A}}(A_2^k)$, $\forall\,k\in\mathbb{N}$.
In particular, it is obvious that abelian approximate self-commutators have
quasitrace zero.
\end{definition}

With this terminology, one has the following result.
\begin{theorem}
Let $\mathcal{A}$ be an AW*-factor of type \twone, and let
$X\in\mathcal{A}_{sa}$ be an element with $q_{\mathcal{A}}(X)=0$.
If
$D\big(\mathbf{s}(X)\big)<1$, then $X$ can be written as
a sum $X=X_1+X_2$, where $X_1,X_2$ are two commuting abelian
approximate self-commutators in $\mathcal{A}$.
\end{theorem}
\begin{proof}
Let $P=I-\mathbf{s}(X)$.
Using the proof of Theorem 5.2 from \cite{KN}, there exist
elements $A_1,A_2,B_1,B_2,Y_1,Y_2,S_1,S_2\in\mathcal{A}_{sa}$, with the
following properties:
\begin{itemize}
\item[(i)] $A_1,A_2,B_1,B_2,Y_1,Y_2,S_1,S_2$ all commute;
\item[(ii)] $A_1\sim B_1$, $A_2\sim B_2$, $Y_1\sim S_1$, $Y_2\sim S_2$, and
$S_1+S_2$ is spectrally symmetric, i.e. $(S_1+S_2)\sim -(S_1+S_2)$;
\item[(iii)] $A_1\perp A_2,B_1,P$, $A_2\perp B_2,P$, $B_1\perp B_2,P$, and
$B_2P=PB_2=Y_1+Y_2$;
\item[(iv)] $Y_1,Y_2\perp S_1,S_2$;
\item[(v)] $Y_1,Y_2,S_1,S_2\in P\mathcal{A}P$;
\item[(vi)] $X=A_1-B_1+A_2-B_2+Y_1+Y_2$.
\end{itemize}
Consider then the elements
\begin{align*}
V_1&=A_1+Y_1-S_2; &V_2&=A_2+\tfrac 12(S_1+S_2);\\
W_1&=B_1+S_1-Y_2; &W_2&=B_2-\tfrac 12(S_1+S_2).
\end{align*}
Using the orthogonal additivity of approximate unitary equivalence
(Corollary 2.1 from \cite{KN}), and the above conditions, it follows that
$V_1\sim W_1$ and $V_2\sim W_2$. Since $V_1,V_2,W_1,W_2$ all commute, it follows that
the elements $X_1=V_1-W_1$ and $X_2=V_2-W_2$ are abelian approximate
self-commutators, and they commute.
Finally, one has
$X_1+X_2=A_1-B_1+A_2-B_2+Y_1+Y_2= X$.
\end{proof}

\begin{corollary}
Let $\mathcal{A}$ be an AW*-factor of type \twone, and let
$X\in\mathcal{A}_{sa}$ be an element with $q_{\mathcal{A}}(X)=0$. There exist
two commuting abelian approximate self-commutators
$X_1,X_2\in \text{\rm Mat}_2(\mathcal{A})$ -- the $2\times 2$ matrix algebra --
such that
\begin{equation}
X_1+X_2=\left[\begin{array}{cc} X &0\\0&0\end{array}\right].
\label{cor2by2}\end{equation}
\end{corollary}
\noindent(According to Berberian's Theorem (see \cite{Be}), the matrix algebra
$\text{\rm Mat}_2(\mathcal{A})$ is an AW*-factor of type \twone.)
\begin{proof}
Denote the matrix algebra $\text{\rm Mat}_2(\mathcal{A})$ by
$\mathcal{A}_2$, and let $\tilde{X}\in \mathcal{A}_2$ denote the matrix in
the right hand side
of \eqref{cor2by2}. It is obvious that, if we consider the projection
$$
E=\left[\begin{array}{cc} 1 &0\\
0&0\end{array}\right],
$$
then $\mathbf{s}(\tilde{X})\leq E$. Since $D_{\mathcal{A}_2}(E)=\frac 12<1$,
and $q_{\mathcal{A}_2}(\tilde{X})=\frac 12 q_{\mathcal{A}}(X)=0$, the desired
conclusion follows imediately from Theorem 2.3
\end{proof}

\end{document}